\input amstex
\documentstyle{amams} 
\input amssym.def
\input amssym.tex

\annalsline{158}{2003}
\received{February 15, 2002}
\startingpage{1041}
\def\A{\mathop{\rm Aut}}

\def\GL{\mathop{\rm GL}}

\def\SL{\mathop{\rm SL}}

\def\End{\mathop{\rm End}}
\def\Sy{\mathop{\rm Sym}}
\def\kb{K}
\def\cd{\cdot }
\def\T{{\rm T}}
\def\St{\mathop{\rm St}}
\def\ch{\mathop{\rm char}}
\def\esd{\mathop{\rm ed}}
\def\s{\frak S_n}

\def\g{\hskip-6pt\mathbold{\gamma}}
\def\de{\hskip-6pt\mathbold{\delta}}

\title{Automorphism groups of \\
finite dimensional simple algebras}  
\shorttitle{Finite dimensional simple algebras}  
\acknowledgements{Both authors were supported in part by
The Erwin Schr\"odinger
International Institute for
Mathematical Physics
(Vienna, Austria).}
 \twoauthors{Nikolai L. Gordeev}{Vladimir L. Popov}
\institutions{ Russian State Pedagogical University,
  St. Petersburg,  Russia\\
{\eightpoint {\it E-mail address\/}: nickgordeev\@mail.ru}\\
\vglue6pt
 Steklov Mathematical Institute,
Russian Academy of Sciences,
Moscow, Russia\\
{\eightpoint {\it E-mail address\/}: popovvl\@orc.ru}\\
\vglue12pt}

\vfill
\centerline{\bf Abstract}
\vglue12pt
 We show that if a field $k$
contains sufficiently many elements
(for instance, if $k$ is infinite), and
$K$ is an algebraically closed field
containing $k$, then every linear
algebraic $k$-group over $K$ is
$k$-isomorphic to $\A(A\otimes_k\!K)$,
where $A$ is a finite dimensional
simple algebra over $k$.
 
\vglue24pt
\section{Introduction}

In this paper, `algebra' over a field means
`nonassociative algebra', i.e., a vector
space $A$ over this field with  
multiplication given by a linear map
$A\otimes A\rightarrow A$, $a_1\otimes
a_2\mapsto a_1a_2$, subject to no {\it a
priori} conditions; cf.\  \cite{Sc}.

Fix a field $k$ and an algebraically closed
field extension $K/k$. Our point of view of
algebraic groups is that of \cite{Bor},
\cite {H}, \cite{Sp}; the underlying
varieties of linear algebraic groups will be
the affine algebraic varieties over $K$. As
usual, algebraic group (resp., subgroup,
homomorphism) defined over $k$ is
abbreviated to $k$-group (resp.,
$k$-subgroup, $k$-homomorphism). If $E/F$ is
a field extension and $V$ is a vector space
over $F$, we denote by $V_E$ the vector
space $V\otimes_F E$  over~$E$.

Let $A$ be a finite dimensional algebra
over $k$ and let $V$ be its underlying
vector space. The $k$-structure $V$ on
$V_K$ defines the $k$-structure on the
linear algebraic group  $\GL(V_K)$. As
$\A (A_{\kb})$, the full automorphism
group of $A_K$, is a closed subgroup of
$\GL(V_K)$, it has the structure of a
linear algebraic group. If $\A
(A_{\kb})$ is defined over $k$ (that is
always the case if $k$ is perfect; cf.\
 [Sp, 12.1.2]), then for each field
extension $F/k$ the full automorphism
group $\A(A_F)$ of $F$-algebra $A_F$ is
the group $\A (A_{K})(F)$ of
$F$-rational points of the algebraic \pagebreak
group $\A (A_{\kb})$.

Let $\Cal A_k$ be the class of linear
algebraic $k$-groups $\A (A_{\kb})$ where
$A$ ranges over all finite dimensional
simple algebras over $k$ such that $\A
(A_{\kb})$ is defined over $k$. It is well
known that many important algebraic groups
belong to $\Cal A_k$: for instance, some
finite simple groups (including the Monster)
and simple algebraic groups appear in this
fashion; cf.\ \cite{Gr}, \cite{KMRT},
\cite{Sp}, \cite{SV}.
 Apart from the
`classical' cases, people studied the
automorphism groups of `exotic' simple
algebras as well; cf.\ \cite{Dix} and
discussion and references in \cite{Pop}. The
new impetus stems from invariant theory: for
$k=\kb$, $\ch k=0$, in \cite{Ilt} it was
proved that if a finite dimensional simple
algebra $A$ over $k$ is generated by $s$
elements, then the field of rational $\A
(A)$-invariant functions of $d\geqslant s$
elements of $A$ is the field of fractions of
the trace algebra (see \cite{Pop} for a
simplified proof). This yields close
approximation to the analogue of classical
invariant theory for some modules of
nonclassical groups belonging to $\Cal A_k$
(for instance, for all simple
$\mathop{\rm E}_8$-modules);
cf.\ \cite{Pop}.

So $\Cal A_k$ is the important class. For
$k=\kb$, $\ch k=0$, it was asked in
\cite{K1} whether all groups in $\Cal A_k$
are reductive. In \cite{Pop} this question
was answered in the negative \footnote{$^1$}{It
was then asked in \cite{K2} whether, for a
simple algebra $A$, the group $\A(A)$ is
reductive if the trace form $(x, y)\mapsto
{\operatorname {tr}}\, L_xL_y$ is
nondegenerate (here and below $L_a$ and
$R_a$ denote the operators of left and right
multiplications of $A$ by $a$). The answer
to this question is negative as well: one
can verify that for some of the simple
algebras with nonreductive automorphism
group constructed in \cite{Pop} all four
trace forms $(x, y)\mapsto {\operatorname
{tr}}\, L_xL_y$, $(x, y)\mapsto
{\operatorname {tr}}\, R_xR_y$, $(x,
y)\mapsto {\operatorname {tr}}\, L_xR_y$ and
$(x, y)\mapsto {\operatorname {tr}}\,
R_xL_y$ are nondegenerate (explicitly, in
the notation of \cite{Pop, (5.18)}, this
holds if and only if $\alpha_1\neq 0$).} 
 and
the general problem of finding a group
theoretical cha\-rac\-te\-rization of $\Cal
A_k$ was raised; in particular it was asked
whether each finite group belongs to $\Cal
A_k$. Notice that each abstract group is
realizable as the full automorphism group of
a (not necessarily finite) field extension
$E/F$, and each finite abstract group is
realizable as the full automorphism group of
a finite (not necessarily Galois) field
extension $F/\Bbb Q$, cf.\ \cite{DG},
\cite{F}, \cite{Ge}.

In this paper we give the complete
solution to the formulated problem. Our
main result is the  following.

\nonumproclaim{ Theorem 1}   If $k$ is a
field containing sufficiently many
elements $($for instance{\rm ,} if $k$ is
infinite{\rm ),} then for each linear
algebraic $k$\/{\rm -}\/group $G$ there is a
finite di\-men\-sional simple algebra
$A$ over~$k$ such that the algebraic
group $\A(A_{\kb})$ is defined over $k$
and $k$\/{\rm -}\/isomorphic to $G$. \endproclaim

The constructions used in the proof of
Theorem 1 yield a precise numerical form of
the condition `sufficiently many'. Moreover,
actually we show that the algebra $A$ in
Theorem 1 can be chosen absolutely simple
(i.e., $A_F$ is simple for each field
extension $F/k$).

From Theorem 1 one immediately deduces
the following corollaries.

{\elevensc Corollary  1.}{\it  Under the
same condition on $k${\rm ,} for each linear
algebraic $k$\/{\rm -}\/group $G$ there is  a
finite dimensional simple algebra $A$
over $k$ such that $G(F)$ is isomorphic
to $\A(A_F)$ for each field extension
$F/k$.  }

\nonumproclaim{{C}orollary  2}  Let $G$ be
a finite abstract group.  Under the
same condition on $k${\rm ,} there is a
finite dimensional simple algebra $A$
over~$k$ such that $\A (A_F)$ is
isomorphic to $G$ for each field
extension $F/k$. 
\endproclaim

One can show (see Section~7) that each
linear algebraic $k$-group can be
realized as the stabilizer of a
$k$-rational element of an algebraic
$\GL(V_{\kb})$-module $M$ defined over
$k$ for some finite dimensional vector
space $V$ over $k$. Theorem~1 implies
that such $M$ can be found among
modules of the very special type:

\nonumproclaim{Theorem 2}  If $k$ is a
field containing sufficiently many
elements $($for instance{\rm ,} if $k$ is
infinite{\rm ),} then for each linear
algebraic $k$\/{\rm -}\/group $G$ there is  a
finite dimensional vector space $V$
over $k$ such that the
$\GL(V_K)$\/{\rm -}\/stabilizer of some
$k$\/{\rm -}\/rational tensor in $V^*_K\!\otimes
\!V^*_{K}\!\otimes \!V^{}_{K}$ is
defined over $k$ and $k$\/{\rm -}\/isomorphic to~$G$. 
\endproclaim

Regarding Theorem 2 it is worthwhile to
notice that $\GL(V_{\kb})$-stabilizers
of points of some dense open subset of
$V^*_{\kb}\!\otimes
\!V^*_{\kb}\!\otimes \!V^{}_{\kb}$ are
trivial; cf.\  \cite{Pop}.

Another application pertains to the notion
of essential dimension. Let $k=\kb$ and let
$A$ be a finite dimensional algebra over
$k$. If $F$ is a field of algebraic
functions over $k$, and $A'$ is an
$F/k$-form of $A$ (i.e., $A'$ is an algebra
over $F$ such that for some field extension
$E/F$ the algebras $A_E$ and $A'_E$ are
isomorphic), put
$$
\zeta(A'):=\min_{F_0}\{{\rm trdeg}_k F_0
\mid A' \text{ is defined over
 the subfield $F_0$
of $F$ containing $k$}\}.
 $$

Define the {\it essential dimension} $\esd
A$ {\it of algebra} $A$ by
$$
\esd A:=\max_F
\max_{A'} \zeta(A').
$$

On the other hand, there is the notion of essential
dimension for algebraic groups  
introduced and studied (for $\ch k=0$) in [Re]. The results in [Re] show that the
essential dimension of $\A(A_{\kb})$
coincides with $\esd A$ and 
demonstrate  how this fact can be used for
finding bounds of essential dimensions of
some linear algebraic groups. The other side of this topic is that
  many (Galois) cohomological invariants
of algebraic groups are defined via
realizations of groups as the automorphism
groups of some finite dimensional algebras,
cf.\ \cite{Se}, \cite{KMRT}, \cite{SV}.
These invariants are the means for finding
bounds of essential dimensions of algebraic
groups as well; cf.\ \cite{Re}. Theorem 1
implies that the essential dimension of {\it
each} linear algebraic group is equal to
$\esd A$ for some simple algebra $A$ over
$k$.

Also, by [Se, Ch.\,III,\,1.1],
Theorem 1 reduces finding Galois
cohomology of {\it each} algebraic
group to describing forms of the
corresponding simple algebra.

Finally, there is the application of our
results to invariant theory as explained
above. For the normalizers $G$ of linear
subspaces in some modules of unimodular
groups our proof of Theorem 1 is
constructive, i.e., we explicitly construct
the corresponding simple algebra $A$ (we
show that every algebraic group is
rea\-lizable as such a normalizer but our
proof of this fact is not constructive).
Therefore for such $G$ our proof yields
constructive description of some $G$-modules
that admit the close approximation to the
analogue of classical invariant theory (in
particular they admit constructive
description of generators of the field of
rational $G$-invariant functions). However
for Corollary~2 of Theorem 1 we are able to
give another, constructive proof (see
Section~5).

Given all this we hope that our results may
be the impetus to finding new interesting
algebras, cohomological invariants, bounds
for essential dimension, and modules that
admit the close approximation to the
analogue of classical invariant theory.

The paper is organized as follows. In
Section 2 for a finite dimensional
vector space $U$ over $k$, we construct
(assuming that $k$ contains
sufficiently many elements) some
algebras whose full automorphism groups
are $\SL(U)$-normalizers of certain
linear subspaces in the tensor algebra
of $U$. In Section~3 we show that the
group of $k$-rational points of each
linear algebraic $k$-group appears as
such a normalizer. In Section 4 for
each finite dimensional algebra over
$k$, we construct a finite dimensional
simple algebra over $k$ with the same
full automorphism group. In Section 5
the proofs of Theorems 1, 2 are given.
In Section 6 we give the constructive
proof of Corollary 2 of Theorem 1.
Since the topic of realizability of
groups as stabilizers and normalizers
is crucial for this paper, for the sake
of completeness we prove in the appendix
(Section~7) several additional results
in this direction.

In September 2001 the second author
delivered a talk on the results of this
paper  at The Erwin Schr\"odinger
International Institute (Vienna).

\demo{Acknowledgement} We are grateful to W.
van der Kallen for useful correspondence and
to the referee  for remarks.
\enddemo

{\it Notation{\rm ,} terminology and conventions}.

\phantom{duh}
\vglue-8pt 

$\bullet $ \ $|X|$ is the number
of elements in a finite set $X$.

\phantom{duh}
\vglue-8pt 

$\bullet $ \  $\A (A)$ is the full
automorphism group of an algebra $A$.

\phantom{duh}
\vglue-8pt 

$\bullet$ \ $\mathop{\rm vect}(A)$ is
the underlying vector space of an
algebra $A$.  

\phantom{duh}
\vglue-8pt 

$\bullet $ \ $K[X]$ is the algebra of
regular function of an algebraic
variety $X$.
\vglue4pt

$\bullet $ \ $\s $ is the symmetric
group of the set $\{1, {\ldots}\,,
n\}$.

$\bullet $ \ $\langle S \rangle $ is
the linear span of a subset $S$ of a
vector space.

\phantom{duh}
\vglue-8pt

$\bullet $ \   Let $V_i$, $i\in I$, be the
vector spaces over a field. When we consider
$V_j$ as the linear subspace of
$\oplus_{i\in I} V_i$, we mean that $V_j$ is
replaced with its copy given by the natural
embedding $V_j\hookrightarrow \oplus_{i\in
I}V_i$. We denote this copy also by $V_j$ in
order to avoid bulky notation; as the
meaning is always clear from the contents,
this does not lead to confusion.

\phantom{duh}
\vglue-8pt

$\bullet $ \ For a finite dimensional
vector space $V$ over a field $F$ we
denote by $\T(V)$ (resp.\ $\Sy(V)$) the
tensor (resp.\ symmetric) algebra of
$V$, and by $\T(V)_+$ (resp.\ $\Sy(V)_+$) its maximal homogeneous
ideal with respect to the natural
grading, $$ \textstyle \T(V)_+:=
\bigoplus_{i\geqslant 1}V^{\otimes i},
\quad \Sy(V)_+:= \bigoplus_{i\geqslant
1}\Sy^i\,(V), \tag1.1 $$
 endowed with the natural
$\GL(V)$-module structure: $$
\textstyle g\cd t_i\!:=\! g^{\otimes
i}(t_i),\ g\cd s_i\!:=\!
\Sy^i(g)(s_i),\ \
 \ g\!\in\! \GL(V),\
t_i \!\in\! V^{\otimes i}\!, \
s_i\!\in\! \Sy^i(V). \tag1.2 $$
  The $\GL(V)$-actions on
$\T(V)$ and $\Sy(V)$ defined by (1.2)
are the faithful actions by algebra
automorphisms. Therefore we may (and
shall) identify $\GL(V)$ with the
corresponding subgroups of $\A(\T(V))$
and $\A(\Sy(V))$.

\phantom{duh}
\vglue-8pt 

$\bullet $ \ For  the finite dimensional vector spaces $V$ and $W$ over a field $F$, a nondegenerate bilinear
pairing $\Delta: V\times W\rightarrow F$ and
a linear operator $g\in \End(V)$ we denote
by $g^*\in \End(W)$ the conjugate of $g$
with respect to $\Delta$.

\phantom{duh}
\vglue-8pt

$\bullet$ \ $\Delta_E$ denotes the
bilinear pairing obtained from
$\Delta$ by a field extension $E/F$.

\phantom{duh}
\vglue-8pt

$\bullet$ \ For a linear operator $t\in
\End(V) $ the eigenspace of $t$
corresponding to the eigenvalue $\alpha $ is
the nonzero linear subspace $\{v\in V \mid
t(v)=\alpha t \}$.

\phantom{duh}
\vglue-8pt 

$\bullet $ \ If a group $G$ acts on a
set $X$, and $S$ is a subset of $X$, we
put $$ G\!_S:=\{g\in G \mid g(S)=S\};
\tag1.3 $$
this is a subgroup of $G$
called the {\it normalizer of} $S$ {\it
in} $G$.

\phantom{duh}
\vglue-8pt 

$\bullet $ \ `Ideal' means `two-sided
ideal'. `Simple algebra' means algebra
with a nonzero multiplication and
without proper ideals. `Algebraic
group' means `linear algebraic group'.
`Module' means `algebraic (`rational'
in terminology of \cite{H}, \cite{Sp})
module'.

\section{Some special algebras}

Let $F$ be a field. In this section we
define and study some finite dimensional
algebras over $F$ to  be used in the
proof of our main result.

\demo{Algebra $A(V, S)$}
Let $V$ be a nonzero finite dimensional
vector space\break over $F$. Fix an integer $r>1$.
Let $S$ be a linear subspace of $V^{\otimes
r}$, resp.\ $\Sy^r(V)$. Then $$
 I(S):=\cases
S\oplus(\bigoplus_{i> r} V^{\otimes i})&
\text{ if } S\subseteq V^{\otimes r},\\
S\oplus(\bigoplus_{i> r} \Sy^{i}(V))& \text{
if } S\subseteq \Sy^{r}(V)
\endcases
\tag2.1 $$
is the ideal of $\T(V)_+$,
resp.\ $\Sy(V)_+$. By definition, $A(V,
S)$ is the factor algebra modulo this
ideal, $$ \gathered A(V,S):=
{\rm A}_+/I(S), \ \text{ where
} {\rm A}_+\!:=\! \cases
\T(V)_+& \text{ if } S\subseteq
V^{\otimes r},\\ \Sy(V)_+& \text{ if }
S\subseteq \Sy^r(V).
\endcases
\endgathered \tag2.2 $$ 
 It readily follows from the
definition that $A(V, S)_E=A(V_E, S_E)$
for each field extension $E/F$.

By (1.1), (2.1), there is natural
identification of graded vector spaces
$$ \textstyle \gathered
\mathop{\rm vect}(A(V, S))=\cases
\hskip -1mm
(\bigoplus_{i=1}^{r-1}V^{\otimes i})
\oplus (V^{\otimes r}/S)& \text{ if }
S\subseteq V^{\otimes r},\\\noalign{\vskip4pt} \hskip -1mm
(\bigoplus_{i=1}^{r-1}\Sy^{i}(V))
\oplus (\Sy^{r}(V)/S)& \text{ if }
S\subseteq \Sy^{r}(V).
\endcases
\endgathered \tag2.3 $$

Restriction of action (1.2) to $\GL(V)_S$
yields a $\GL(V)_S$-action on
${\rm A}_+$. By (2.1), the ideal
$I(S)$ is $\GL(V)\!_S$-stable. Hence (2.2)
defines a $\GL(V)_S$-action on $A(V, S)$ by
algebra automorphisms, and the canonical
projection $\pi\!$ of ${\rm A}_+$
to $A(V, S)$ is $\GL(V)_S$-equivariant. The
condition $r>1$ implies that $V=V^{\otimes
1}=\Sy^1(V)$ is a submodule of the
$\GL(V)_S$-module $A(V,S)$. Hence
$\GL(V)\!_S$ acts on $A(V, S)$ faithfully,
and we may (and shall) identify $\GL(V)\!_S$
with the subgroup of $\A (A(V, S))$.
\enddemo

\nonumproclaim{Proposition 1}  
$\{\sigma\in \A (A(V, S)) \mid \sigma
(V)=V\}=\GL(V)\!_S$. 
\endproclaim

\demo{Proof} It readily follows from
(1.1)--(2.2)  that the right-hand side
of this equa\-li\-ty is contained in
its left-hand side.

To prove the inverse inclusion, take an
element $\sigma\in \A (A(V, S))$ such
that $\sigma (V)= V$. Put $g:=\sigma
|_V$. Consider $g$ as the automorphism
of ${\rm A}_+$ defined by
(1.2). We claim that the diagram $$ \CD
{\rm A}_+ @>g>>
{\rm A}_+
\\ @V{\pi}VV @VV{\pi}V\\ A(V,
S)@>{\sigma}>>A(V, S)
\endCD,
\tag2.4 $$

\noindent cf.\ (2.2), is commutative.
To prove this, notice  that as the
algebra ${\rm A}_+$ is
gene\-ra\-ted by its homogeneous
subspace $V$ of degree 1 (see (2.3)),
it suffices to check the equality
$\sigma (\pi(x))=\pi(g(x))$ only for
$x\in V$. But in this case it is
evident since $g(x)=\sigma(x)\in V$ and
$\pi(y)=y$ for each $y\in V$.

Commutativity of (2.4) implies that $g
\cdot \ker \pi = \ker \pi$. As  $\ker
\pi = I(S)$, for\-mu\-las (2.1), (1.2),
(1.3) imply that $g\in \GL(V)\!_S$.
Hence $g$ can be considered as the
automorphism of $A(V,S)$ defined by
(2.2). Since its restriction to the
subspace $V$ of $A(V,S)$ coincides with
that of $\sigma$, and $V$ generates the
algebra $A(V,S)$, this automorphism
coincides with $\sigma$, whence
$\sigma\in \GL(V)_S$.  
\enddemo

 {\it Algebra $B(U)$.}
Let $U$ be a nonzero finite dimensional
vector space over~$F$, and $n:=\dim U$.

The algebra $B(U)$ over $F$ is defined
as follows. Its underlying vector space
is that of the exterior algebra of $U$,
$$ \textstyle \mathop{\rm vect}(B(U))
=\bigoplus_{i=1}^{n} \wedge^i U.
\tag2.5 $$

To define the multiplication in $B(U)$
fix a basis of each $\wedge^i U$,
$i=1,{\ldots\,},n$. For $i=n$, it
consists of a single element $b_0$. The
(numbered) union of these bases is a
basis $\Cal B_{B(U)}$ of
$\mathop{\rm vect}(B(U))$. By
definition, the multiplication in
$B(U)$ is given by $$ pq=\cases p\wedge
q, & \text{ for \hskip 0,7mm $p, q \in
{\Cal B}_{B(U)}$,} \text{ and \hskip
0,7mm $p$ or $q\neq b_0$},\\ \hskip
2,8mm b_0, & \text{ for \hskip 0,7mm
$p=q=b_0$.}
\endcases
\tag2.6 $$
It is immediately seen that
up to isomorphism $B(U)$ does not
depend on the choice of $\Cal
B_{B(U)}$, and $ B(U)_E=B(U_E)$ for
each field extension $E/F$.

The $\GL(U)$-module structure on $\T(U)$
given by (1.2) for $V=U$ induces the
$\GL(U)$-module structure on
$\mathop{\rm vect}(B(U))$ given by $$ g\cd
x_i= (\wedge^i g)(x_i), \hskip 2,5mm g\in
\GL(U),\ x_i\in \wedge^i U. \tag2.7 $$
In particular $$ g\cd b_0=(\det
g)b_0, \hskip 2,5mm g\in \GL(U). \tag2.8 $$
As $U=\wedge^1U$ is the submodule of
$\mathop{\rm vect}(B(U))$, the
$\GL(U)$-action on
$\mathop{\rm vect}(B(U))$ is faithful.
Therefore we may (and shall) identify
$\GL(U)$ with the subgroup of
$\GL(\mathop{\rm vect}(B(U)))$.

\nonumproclaim{Proposition 2}  
$\{\sigma \in \A (B(U)) \mid \sigma
(U)=U\}=\SL(U)$. 
\endproclaim

\demo{Proof} First we show that the
left-hand side of this equality is
contained in its right-hand side. Take
an element $\sigma\in\A (B(U))$ such
that $\sigma(U)=U$. By (2.5) and (2.6),
the algebra $B(U)$ is generated by $U$.
Together with (2.5), (2.6), (2.7), this
shows that $\sigma(x)=\sigma|_U\cd x$
for each element $x\in B(U)$. In
particular, $\sigma (b_0)\in \wedge^n
U$. As  $\sigma$ is an automorphism of
the algebra $B(U)$, it follows from
(2.6) that $b_0$ and $\sigma(b_0)\in
\wedge^n U$ are the idempotents of this
algebra. But $\dim \wedge^n U=1$
readily implies that $b_0$ is the
unique idempotent in $\wedge^n U$.
Hence $\sigma|_U\cd b_0=b_0$. By (2.8),
this gives $\sigma\in \SL(U)$.

Next we show that the right-hand side
of the equality under the proof is
contained in its left-hand side. Take
an element $g\in \SL(U)$ and the
elements $p, q\in {\Cal B}_{B(U)}$. We
have to prove that $$ g\cd (pq)
=
(g\cd p)(g\cd q). \tag2.9 $$

\noindent Let, say, $p\neq b_0$. Then
by (2.7) we have $g\cd p=\sum_{b\in
{\Cal B}_{B(U)}}\alpha_{b}b$ for some
$\alpha_b\in F$ where $\alpha_{b_0}=0$.
By (2.6) and (2.7), we have $g\cd
(pq)=g\cd(p\wedge q)= (g\cd p)\wedge
(g\cd q)= \sum_{b\in {\Cal
B}_{B(U)}}\alpha_{b}(b\wedge (g\cd q))=
\sum_{b\in {\Cal
B}_{B(U)}}\alpha_{b}(b(g\cd q))=
(\sum_{b\in {\Cal
B}_{B(U)}}\alpha_{b}b)(g\cd q) =(g\cd
p)(g\cd q).$ Thus (2.9) holds in this
case. Then similar arguments show that
(2.9) holds for $q\neq b_0$. Finally,
from (2.6), (2.8) and $\det g=1$ we
obtain $g\cd (b_0b_0)=g\cd b_0 = b_0 =
b_0b_0= (g\cd b_0)(g\cd b_0)$. Thus
(2.9) holds for $p=q=b_0$ as well.
\enddemo

{\it Algebra $C(L,U,\g)$.}

\nonumproclaim{Lemma 1}  Let $A$ be an algebra
over $F$ with the left identity $e\in A$
such that $\mathop{\rm vect}(A)=\langle e
\rangle \oplus A_1\oplus{\cdots} \oplus
A_r,$ where $A_i$ is the eigenspace with a
nonzero eigen\-value $\alpha_i$ of the
operator of right multiplication of $A$ by
$e$. Then \medbreak \item{\rm (i)} $e$ is
the unique left identity in $A${\rm ;} 
\smallbreak \item{\rm
(ii)} if $\sigma\in \A (A)${\rm ,} then $\sigma
(e)=e$ and $\sigma(A_i)=A_i$ for all $i$.
\endproclaim

\demo{Proof} 
(i) Let $e'$ be a left identity of $A$. As
$e'=\alpha e + a_1+{\cdots}+a_r$ for some
$\alpha\in F$, $a_i\in A_i$, we have $e=
e'e=(\alpha e + a_1+{\cdots}+a_r)e= \alpha e
+ \alpha_1a_1+{\cdots}+\alpha_ra_r$. Since
$\alpha_i\neq 0$ for all $i$, this implies
$\alpha=1$ and $a_i=0$ for all $i$, i.e.,
$e'=e$.
\vglue4pt
(ii) As $\sigma(A_i)$ is the eigenspace with
the eigenvalue $\alpha_i$ of the operator of
right multiplication of $A$ by $\sigma(e)$,
and $1\neq\alpha_i\neq\alpha_j$ for all $i$
and $j\neq i$ because of the definition of
eigenspace (cf.\;Intro\-duc\-tion), (ii)
follows from (i).  
\enddemo

Fix two nonzero finite dimensional vector
spaces $L$ and $U$ over $F$. Put $s:=\dim
L$, $n:=\dim U$ and assume that
$$|F|\geqslant \max\{n+3, s+1\}.\tag2.10$$

\nonumproclaim{Lemma 2}  There is a
structure of $F$\/{\rm -}\/algebra on $L$ such
that $\A
(L_E)=\{{\rm id}_{L_E}\}$ for
each field extension $E/F$. 
\endproclaim

\demo{Proof} If $s=1$, each nonzero
multiplication on $L$ gives the structure we
are after. If $s>1$, consider a basis $e,
e_1, \ldots, e_{s-1}$ of $L$ and fix any
algebra structure on $L$ satisfying the
following conditions (by (2.10), this is
possible):
\medbreak
 \item{(L1)} $e$ is the
left identity; 
\vglue4pt \item{(L2)} each $\langle
e_i\rangle$ is the eigenspace with a nonzero
eigenvalue of the operator of right
multiplication of $L$ by $e$. 
\vglue4pt\item{(L3)}
$e_i^2\in \langle e_i\rangle \setminus
\{0\}$ for each $i$.
\vglue4pt
\indent By Lemma 1, if $\sigma\in \A (L_E)$,
we have $\sigma(e)=e$ and $\sigma (\langle
e_i\rangle_E)=\langle e_i\rangle_E$ for each
$i$. Whence $\sigma =
{\rm id}_{L_E}$.  
\enddemo

Fix a sequence $\g=
(\gamma_1,\ldots,\gamma_{n+1})\in F^{n+1},
\gamma_i\neq 0, 1, \ \gamma_i\neq \gamma_j
\text { for } i\neq j$; by (2.10), this is
possible. Using Lemma 2, fix a structure of
$F$-algebra on $L$ such that $\A
(L_E)=\{{\rm id}_{L_E}\}$ for each
field extension $E/F$. We use the same
notation $L$ for this algebra.

The algebra $C(L, U, \g)$ is defined as
follows. By definition, the direct sum of
algebras $L$ and $B(U)$ is the subalgebra of
$C(L, U, \g)$, and there is an element $c\in
C(L, U, \g)$ such that \vskip -5mm
$$ \textstyle \mathop{\rm vect}(C(L,
U, \g))=\langle c\rangle\oplus
\mathop{\rm vect}(L\oplus B(U))
\overset\text{(2.5)}\to= \langle
c\rangle\oplus
\mathop{\rm vect}(L)\oplus
(\bigoplus_{i=1}^{n}\wedge^i U)
\tag2.11 $$

\noindent and the following conditions hold:
\medbreak \item{(C1)} $c$ is the left identity
of $C(L, U, \g)$;
\item{(C2)}$\mathop{\rm vect}(L)$ and
$\wedge^i U$, $i=1,{\ldots\,}, n$, in (2.11)
are respectively the eigenspaces with
eigenvalues $\gamma_1,\dots,\gamma_{n+1}$ of
the operator of right mutiplication of $C(L,
U, \g)$ by $c$.
\vglue4pt

It is immediately seen that $C(L, U,
\g)_E=C(L_E, U_E, \g)$ for each field
extension $E/F$.

Define the $\GL(U)$-module structure on
$\mathop{\rm vect}(C(L, U, \g))$ by the
condition that in (2.11) the subspaces
$\langle c\rangle$ and
$\mathop{\rm vect}(L)$ are trivial
$\GL(U)$-submodules, and
$\oplus_{i=1}^{n}\wedge^i U$ is the
$\GL(U)$-submodule with $\GL(U)$-module
structure defined by (2.7). Thus for all
$g\in \GL(U),\ \alpha\in F,\ l\in L,\ x_i\in
\wedge^i U$, $$ \textstyle g\cd (\alpha c +
l +\sum_{i=1}^n x_i) = \alpha c + l +
\sum_{i=1}^{n}(\wedge^i g)(x_i). \tag2.12 $$

The $\GL(U)$-action on
$\mathop{\rm vect}(C(L, U, \g))$ given by
(2.12) is faithful.\break Therefore we may (and
shall) consider $\GL(U)$ as the subgroup of\break
$\GL(\mathop{\rm vect}(C(L, U, \g)))$.

\nonumproclaim{Proposition 3}  
$\A (C(L, U, \g ))= \SL(U)$. 
\endproclaim

\demo{Proof} The claim follows from the
next two:

\smallskip

(i) $\A(C(L, U, \g))\subset \GL(U)$;
\vglue4pt
(ii) $g\in \GL(U)$ lies in $\A(C(L, U, \g))$
if and only if $g\in \SL(U)$.

\pagebreak

To prove (i), take an element
$\sigma\in \A (C(L, U, \g))$. By
(2.12), we have to show that all direct
summands in the right-hand side of
$(2.11)$ are $\sigma$-stable, and
$\sigma(x)=\sigma|_{U}\cd x$ for all
$x\in C(L, U, \g)$. The first statement
follows from (C1), (C2) and Lemma~1. As
$L$ and $B(U)$ are the subalgebras of
$C(L, U, \g)$, the second statement
follows from the first, the condition
$\A (L)=\{{\rm id}_L\}$, Lemma
1, Proposition 2 and formulas (2.7),
(2.12).

To prove the `only if' part of (ii), assume
that $ g\in \A (C(L, U, \g))$. As by (i) and
(2.12) the subalgebra $B(U)$ is $g$-stable,
$g|_{B(U)}$ is a well-defined element of  $\A
(B(U))$. By Proposition 2, this gives $g\in
\SL(U)$.

To prove the `if' part of (ii) assume
that $g\in \SL(U)$. By (2.12) the
subalgebra $L\oplus B(U)$ of $C(L, U,
\g )$ is $g$-stable, and by Proposition
2 the transformation $g|_{L\oplus
B(U)}$ is its automorphism. Hence it
remains to show that if $x$ is an
element of some direct summand of the
right-hand side of (2.11), then the
following equalities hold: $$ g\cd (cx)
= (g\cd c)(g\cd x) \quad {\text {\rm
and}} \quad g\cd (xc) = (g\cd x)(g\cd
c).$$

By (2.12), $ g\cd c =c$. Together with
(C1) this yields the first equality needed: $ g\cd (cx)= g\cd x =c (
g\cd x)= (g\cd c) ( g\cd x)$. By (C2),
$xc=\alpha x$ for some $\alpha\in F$.
Hence (C1), (C2) and (2.12) imply
$
(g\cd x)c=\alpha ( g\cd x)$. This
yields the second equality needed:
$
g\cd (xc)= g\cd (\alpha x)=\alpha (g\cd
x)= (g\cd x)c= (g\cd x) (g\cd c)$.
\enddemo

{\it Algebra $D(L, U, S, \g, \de,
\Phi)$.}
Let $L$ and $U$ be two nonzero
finite dimensional vector spaces over
$F$. Put $s:=\dim L$, $n:=\dim U$ and
$$ V:=L\oplus U. \tag2.13 $$

\noindent Let $r>1$ be an integer.
Assume that $$|F|\geqslant \max\{n + 3,
s + 1, r + 3\}\tag2.14$$

\noindent and fix the following data:

\medbreak \item{(i)} a linear subspace $S$ of
$V^{\otimes r}$; 
\vglue6pt\item{(ii)} two sequences
$\g= (\gamma_1,\ldots ,\gamma_{n+1})\in
F^{n+1}$,
$\de=(\delta_1,\ldots,\delta_{m})\in F^{m}$,
where $m=r+1$ if $S\neq V^{\otimes r}$ and
$m=r$ otherwise, and $\gamma_i, \delta_j\in
F\setminus\{0, 1\}$, $\gamma_i\neq\gamma_j$,
$\delta_i\neq \delta_j$ for  $i\neq j$ (by
(2.14), this is possible);
\vglue6pt
 \item{(iii)} a
structure of $F$-algebra on $L$ such that
$\A (L_E)=\{{\rm id}_{L_E}\}$ for
each field extension $E/F$ (by (2.14) and
Lemma 2, this is possible); we use the same
notation $L$ for this algebra.
\pagebreak

Define the algebra $D(L, U, S, \g, \de,
\Phi)$ as follows. First, $A(V, S)$,\break $C(L,
U, \g)$ are the subalgebras of $D(L, U, S,
\g, \de, \Phi)$ and the sum of their
underlying vector spaces is direct.  Thus
$\mathop{\rm vect}(D(L, U, S, \g, \de,
\Phi))$ contains two distinguished copies of
$V$: the copy $V^{}_A$ corresponds to the
summand $V^{\otimes 1}$ in (2.3),
and the copy $V^{}_C$ to the summand
$\mathop{\rm vect}(L)\oplus \wedge\!^1 U$
in (2.11).

Denote by $L^{}_A$, $U^{}_A$, resp.\ $L^{}_C$, $U^{}_C$, the copies of resp.\ $L$,
$U$ (see (2.13)) in $V^{}_A$,
resp.~$V^{}_C$, and fix a nondegenerate
bilinear pairing $ \Phi:\ V^{}_A\times
V^{}_C\rightarrow F $ such that $L^{}_A$ is
orthogonal to $U^{}_C$, and $U^{}_A$ to
$L^{}_C$,
$$
\Phi|^{}_{L^{}_A\times U^{}_C}=0, \quad
\Phi|^{}_{U^{}_A\times L^{}_C}=0. \tag2.15
$$

Second, there is an element $d\in D(L, U, S,
\g, \de, \Phi)$ such that $$
\mathop{\rm vect}(D(L, U, S, \g, \de,
\Phi))=\langle d \rangle \oplus
\mathop{\rm vect}(A(V, S))\oplus
\mathop{\rm vect}(C(L, U, \g)) \tag2.16 $$
 and the following conditions hold:
\medbreak \item{(D1)} $d$ is the left identity
of $D(L, U, S, \g, \de, \Phi)$;
\vglue6pt \item{(D2)}
$\mathop{\rm vect}(C(L, U, \g))$ in (2.16)
and all summands in decomposition (2.3) of
the summand $\mathop{\rm vect}(A(V, S))$
in (2.16) are the eigenspaces with
eigenva\-lues $\delta_1,\ldots, \delta_{m}$
of the operator of right multiplication of
the algebra\break $D(L, U, S, \g, \de, \Phi)$ by
$d$; 
\vglue6pt\item{(D3)} if $x\in
\mathop{\rm vect}(A(V, S))$ and $y\in
\mathop{\rm vect}(C(L, U, \g))$ are the
elements of some direct summands in the
right-hand sides of (2.3) and (2.11)
respectively, then
 their product in $D(L, U, S, \g,
\de, \Phi)$ is given by $$ xy=yx=\cases
\Phi(x,y)d & \text{ if $x\in V^{}_A,\,
y\in V^{}_C$},\\
\hskip 4mm 0 & \text{ otherwise. }
\endcases
\tag2.17 $$
\vglue4pt

It is immediately seen that for each
field extension $E/F$ $$ D(L, U, S, \g,
\de, \Phi )_E=D(L_E, U_E, S_E, \g, \de,
\Phi_E). \tag2.18$$

We identify $g\in \GL(U)$ with
${\rm id}_L\oplus g\in \GL(V)$ and
consider $\GL(U)$ as the subgroup of
$\GL(V)$. Formulas (2.2), (1.2) , (2.12) and
(1.3) define the $\GL(U)\!_S$-module
structures on $\mathop{\rm vect}(A(V, S))$
and $\mathop{\rm vect}(C(L, U, \g ))$.

\nonumproclaim{Proposition 4}  With
respect to decomposition $(2.16)$ and
the bilinear pairing $\Phi${\rm ,}  $$ \A (D(L,
U, S, \g, \de, \Phi))=
\{{\rm id}_{\langle d
\rangle}\oplus \ g \oplus (g^*)^{-1}
\mid g\in \SL(U)\!_S\}. \tag2.19 $$
\endproclaim \pagebreak

\demo{Proof} Take an element $\sigma \in \A
(D(L, U, S, \g, \de, \Phi))$. Lemma 1 and
conditions (D1), (D2) imply that
$\sigma(d)=d$ and that the summands in
(2.16) and (2.3) are $\sigma$-stable. From
condition (D1) and Propositions 1, 3 we
deduce that $\sigma
={\rm id}_{\langle d \rangle}
\oplus g \oplus h$ for some $g\in
\GL(V)\!_S$, $h\in \SL(U)$.

Let $x\in V^{}_A$, $y \in V^{}_C$. Then
$\sigma(x)=g\cdot x\in V^{}_A$,
$\sigma(y)=h\cdot y\in V^{}_C$. So from
(2.17) we obtain $\Phi(g\cdot x, h\cdot y)d=
\sigma(x)\sigma(y)= \sigma(xy)=
\sigma(\Phi(x, y)d)= \Phi(x, y)\sigma(d)=
\Phi(x, y)d$. Hence $$ \Phi(g\cd x,h\cd y)=
\Phi(x, y),\ \ x\in V^{}_A,\, y\in V^{}_C.
\tag2.20
$$

As $h\cdot L=L$ and $h\cdot U=U$, it follows
from (2.15), (2.20) and nondegene\-ra\-cy of
$\Phi $ that $g\cd L=L$ and $g\cd U=U$.
Besides, (2.15) and nondegeneracy of $\Phi $
show that the pairings $\Phi|_{L^{}_A\times
L^{}_C}$ and $\Phi |_{U^{}_A\times U^{}_C}$
are nondegenerate. Hence by (2.20) $$
h=(g^*)^{-1}, \ g|^{}_L=((h|_L)^*)^{-1}, \
g|_U=((h|^{}_U)^*)^{-1}. \tag2.21 $$

\noindent As $h\in \SL(U)$, we have
$h|^{}_L={\rm id}^{}_L$ and
$\det(h|^{}_U)=1$. By (2.21), this gives
$g|^{}_L={\rm id}_L$ and
$\det(g|^{}_U)=1$, i.e., $g\in \SL(U)$. Thus
$g\in \SL(U)_S$, i.e., the left-hand side of
equality (2.19) is contained in its
right-hand side.

To prove the inverse inclusion, take an
element
$\varepsilon={\rm id}_{\langle d
\rangle}\! \oplus\,q\oplus\! (q^*)^{-1}$,
whe\-re $q\in \SL(U)_S$. We have to show
that $$
\varepsilon(xy)=\varepsilon(x)\varepsilon(y),\
\ x, y \in D(L, U, S, \g, \de, \Phi).
\tag2.22 $$

 By Proposition 1, we have
$\varepsilon|_{A(V, S)}\in \A (A(V,
S))$. By (D1), this gives (2.22) for
$x, y\in A(V, S)$.

As above, we have $(q^*)^{-1}\in
\SL(U)$. Hence $\varepsilon|_{C(L, U,
\g)}\in \A (C(L, U, \g))$ by
Proposition 3. By (D1) this gives
(2.22) for $x, y\in C(L, U, \g)$.

If $x=d$, then (2.22) follows from (D1)
and the equality $\varepsilon(d)=d$. If
$y=d$, then (2.22) follows from (D2) as
$\varepsilon(C(L, U, \g))=C(L, U, \g)$
and each summand in (2.3) is
$\varepsilon$-stable.

Further, let $x$ (resp., $y$) be an element
of some direct summand of $A(V, S)$ in
decomposition (2.3), and $y$ (resp., $x$) be  an
element of some direct summand of $C(L, U,
\g)$ in the right-hand side of (2.16). As
these summands are $\varepsilon$-stable,
(2.11) implies that $xy=yx=0$ and
$\varepsilon(x)\varepsilon(y)=
\varepsilon(y)\varepsilon(x)=0$, and hence
(2.22) holds, unless $x\in V^{}_A$ and $y\in
V^{}_C$ (resp., $x\in V^{}_C$ and $y\in
V^{}_A$).

Finally, let $x\in V^{}_A$, $y\in V^{}_C$.
Then $\varepsilon(x)=q\cdot x$,
$\varepsilon(y)=(q^*)^{-1}\cdot y$, so by
(2.17) we have
$\varepsilon(xy)=\varepsilon(\Phi(x, y)d)
=\Phi(x, y)\varepsilon(d)= \Phi(x ,
y)d=\Phi(x, q^*(q^*)^{-1}\cdot y)d=
\Phi(q\cdot x, (q^*)^{-1}\cdot y)d =
\varepsilon(x)\varepsilon(y)$. Hence (2.22)
holds in this case. This and (2.17) show 
that (2.22) holds for $x\in V^{}_C$, $y\in
V^{}_A$ as well.   
\enddemo

\nonumproclaim{{C}orollary}  \hskip-3pt Assume that $F=K$
and $L${\rm ,} $U${\rm ,} $S${\rm ,} $\Phi$ are defined  over
$k$.\break If $\,\SL(U)_S$ is the $k$\/{\rm -}\/group{\rm ,} then
$\A (D(L, U, S, \g, \de, \Phi))$ is the
$k$\/{\rm -}\/group $k$\/{\rm -}\/iso\-morphic to $\SL(U)_S$. 
\endproclaim

\demo{Proof} As $\A (D(L, U, S, \g,
\de, \Phi))$ is the image of the 
$k$-homo\-mor\-phism of $k$-groups
$\SL(U)_S\rightarrow
\GL(\mathop{\rm vect}(D(L, U, S, \g,
\de, \Phi))), \ g\mapsto
{\rm id}_{\langle d
\rangle}\!\oplus\,g\!\oplus\!
(g^*)^{-1}$, it is the $k$-group as
well; cf.\ \cite{Sp, 2.2.5}. Considered
as the  $k$-homomorphism of $k$-groups
$\SL(U)_S\rightarrow \A (D(L, U, S, \g,
\de, \Phi))$, this $k$-homomorphism is
$k$-isomorphism since there is the
inverse $k$-homomorphism
${\rm id}_{\langle d
\rangle}\!\oplus\,g
\oplus\!(g^*)^{-1}\mapsto g$.  
\enddemo

\vglue-6pt
\section{Normalizers of linear subspaces
in some modules}
\vglue-6pt

In Section 2 we realized normalizers of
linear subspaces in some modules of
unimodular groups as the full automorphism
groups of some algebras. Now we shall show
that each group appears as such a
normalizer.

\nonumproclaim{Proposition 5} Let $G$ be an
algebraic $k$\/{\rm -}\/group. There is a finite
dimensional vector space $U$ over $\kb$
endowed with a $k$\/{\rm -}\/structure{\rm ,} and an integer
$b\geqslant 0$ such that the following
holds. Let $r$ be an integer{\rm ,} $r\geqslant
b${\rm ,} and $L$ be a trivial finite dimensional
$\SL(U)$\/{\rm -}\/module defined over $k${\rm ,}
 $\dim
L\geqslant 2$. Then the $\SL(U)$\/{\rm -}\/module
$(L\oplus U)^{\otimes r}$ contains a linear
subspace $S$ defined over $k$ such that
$\SL(U)_S$ is the $k$\/{\rm -}\/subgroup of $\SL(U)$
$k$-isomorphic to $G$.  
\endproclaim

\demo{Proof} One may realize $G$ as a
closed $k$-subgroup of $\GL_m(K)$ for
some $m$, cf.\  [Sp, 2.3.7], and
obviously $\GL_m(K)$ may be realized as
a closed $k$-subgroup of
$\SL_{m+1}(K)$. Therefore we may (and
shall) consider $G$ as a closed
$k$-subgroup of $\SL(U)$ for some
finite dimensional vector space $U$
over $\kb$ endowed with a
$k$-structure.

Consider the $k$-structure on $\End
(U)$ defined by the $k$-structure of
$U$. Define the $\SL(U)$-module
structure on $\End (U)$ by left
multiplications, $g\!\cdot\!
h\!:=\!g\circ~\hskip -1.5mmh$, $g\in
\SL(U)$, $h\in \End(U)$. The
$\SL(U)$-module $\End (U)$ is defined
over $k$. The subvariety $\SL(U)$ of
$\End(U)$ is closed, defined over $k$
and $\SL(U)$-stable. The restriction to
$\SL(U)$ of the $\SL(U)$-action on
$\End(U)$ is the action by left
translations.

These $\SL(U)$-actions endow the
algebras $K[\SL(U)]$ and $K[\End(U)]$\break
with the structures of al\-geb\-raic
$\SL(U)$-modules defined over $k$.
Restriction of\break functions yields the
$k$-defined $\SL(U)$-equivariant
epimorphism of algebras\break
$\kb[\End(U)]\rightarrow \kb[\SL(U)]$,
$f\mapsto f|_{\SL(U)}$. Since the
$\SL(U)$-module $\End (U)$ is
isomorphic over $k$ to $U^{\oplus d}$,
where $d:=\dim U$, we deduce from here
that there is a $k$-defined
$\SL(U)$-equivariant epimorphism of
algebras $\kb[U^{\oplus d}]\rightarrow
\kb[\SL(U)]$. In turn, as $\kb
[U^{\oplus d}] =\Sy((U^*)^{\oplus d})$,
this and the definition of symmetric
algebra yield that there is a
$k$-defined $\SL(U)$-equivariant
epi\-morphism of algebras $$\alpha:
\T((U^*)^{\oplus d})\rightarrow
\kb[\SL(U)].\tag3.1$$

The classical Chevalley argument,
cf.\  [Sp, 5.5.1], shows that
$\kb[\SL(U)]$ contains a
finite dimensional linear subspace $W$
defined over $k$ such that
$$\SL(U)_W=G.\tag3.2$$

As $W$ is finite dimensional and (3.1)
is an epimorphism, there is an integer $h
>0$ such that $$\textstyle W\subseteq
\alpha (\bigoplus_{i\leqslant
h}((U^{*})^{\oplus d})^{\otimes i}).
\tag3.3$$
Put $W':=\alpha^{-1}(W)\cap
(\bigoplus_{i\leqslant h}((U^{*})^{\oplus
d})^{\otimes i})$. Since $((U^{*})^{\oplus
d})^{\otimes i}$ is $\SL(U)$-stable, (3.3)
implies that
$$\SL(U)_{W'}=\SL(U)_W.\tag3.4$$

Because of $\dim L\geqslant 2$, one can
find in $L$  two linearly independent\break
$k$-rational elements $l_1$ and $l_2$.
We claim that there is a $k$-defined
injection of $\SL(U)$-modules $$\iota:
\T((U^*)^{\oplus d})\hookrightarrow
\T(\langle l_1\rangle \oplus
U^*).\tag3.5$$
To prove this, let $U^*_i$ be the
$i^{\rm th}$ direct summand of  $(U^*)^{\oplus d}$
considered as the linear subspace of
$(U^*)^{\oplus d}$. Fix a basis
$\{f_{ij}\mid j=1,\ldots, d\}$ of $U^*_i$
consisting of $k$-rational elements. For any
$i_1, j_1, \ldots,i_t, j_t\in [1,d]$, $t=1,
2, \ldots $, define the element of
$\T(\langle l_1\rangle \oplus U^*)$ by
$$ \iota(f_{i_1j_1}\otimes
\cdots\otimes f_{i_tj_t}):= l^{\otimes
i_1}_1\otimes f'_{i_1j_1}\otimes \cdots
\otimes l^{\otimes i_t}_1\otimes
f'_{i_tj_t},\tag3.6
$$ where $f'_{ij}$ is the image of $f_{ij}$
under the natural isomorphism
$U^*_i\rightarrow U^*$. Then one easily
verifies that the linear mapping $\iota:
\T((U^*)^{\oplus d})\rightarrow \T(\langle
l_1\rangle\oplus U^*)$ defined on the basis
$\{f_{i_1j_1}\otimes\cdots\otimes
f_{i_tj_t}\}$ of $ \T((U^*)^{\oplus d})_+$
by formula (3.6) and sending $1$ to $1$ has
the properties we are after.

From existence of embedding (3.5) it follows
that retaining the normalizer in $\SL(U)$
one can replace the subspace $W'$ with
another one, $W'':=\iota(W')$:
$$\SL(U)_{W''}=\SL(U)_{W'}.\tag 3.7 $$

Since $W''$ is finite dimensional, there is
an integer $b\geqslant 0$ such that
$$\textstyle W''\subseteq
\bigoplus_{i\leqslant b}(\langle
l_1\rangle\oplus U^*)^{\otimes i}.\tag3.8$$
 Take an integer $r\geqslant
b$ and consider the linear mapping
$$\textstyle \iota_r:
\bigoplus_{i\leqslant b}(\langle
l_1\rangle \oplus U^*)^{\otimes
i}\rightarrow (L\oplus U^*)^{\otimes
r},\ \ f_i\mapsto l_{2}^{\otimes
(r-i)}\otimes f_i,\ f_i\in (\langle
l_1\rangle\oplus U^*)^{\otimes i}.$$
 It is immediately seen that
$\iota_r$ is the injection of
$\SL(U)$-modules defined over $k$. From this
and (3.8) we deduce that the normalizers in
$\SL(U)$  of the subspaces $W''$ and
$\iota_r(W'')$ coincide,
$$\SL(U)_{\iota_r(W'')}=
\SL(U)_{W''}.\tag3.9$$

Now we take into account that there is
a $k$-automorphism $\sigma\in
\A(\SL(U))$ of order $2$ such that the
$\SL(U)$-module $U$ is isomorphic over
$k$ to the\break $\SL(U)$-module with
underlying vector space $U^*$ and
$\SL(U)$-action defined by $$g\ast
f:=\sigma(g)(f), \ \ g\in \SL(U), f\in
U^*\tag 3.10$$
 (fixing a $k$-rational basis
in $U$ and the dual basis in $U^*$,
identify
 $\SL(U)$ with
$\SL_d(K)$, and $U$, $U^*$ with $K^d$; then
$\sigma(g)=(g^{\mathop{\rm T}})^{-1}$).
Hence replacing the standard $\SL(U)$-action
on $(L\oplus U^*)^{\otimes r}$ with the
action defined by (3.10) (and trivial on
$L$) we obtain the $\SL(U)$-module that is
isomorphic to $(L\oplus U)^{\otimes r}$ over
$k$. The normalizer in $\SL(U)$ of the
subspace $\iota_r(W'')$ of this module is
the subgroup
$\sigma(\SL(U)_{\iota_r(W'')})$. Since it is
$k$-isomorphic to $\SL(U)_{\iota_r(W'')}$,
the claim follows from (3.2), (3.4), (3.7)
and (3.9).  
\enddemo

\section{Simple and nonsimple algebras}

Let $R$ be a finite dimensional algebra over
a field $F$. Assume that $|F|\geqslant 4$
(however see the remark at the end of this
section). Fix the following data: \medbreak
\item{(a)} two nonzero elements $\alpha,
\zeta\in F$, $\alpha, \zeta \neq 1$,
$\alpha\neq\zeta$, 
\vglue4pt
\item{(b)} an algebra $Z$
over $F$ of the same dimension as $R$ and
with zero multiplication,
$$
\quad z_1z_2=0 \ \text{ for all } z_1, z_2
\in Z, \tag4.1 $$
\vglue4pt \item{(c)} a nondegenerate
bilinear pairing $\Delta: Z\times
R\rightarrow F$. \vglue4pt

 In this section
we construct a finite dimensional algebra
$R(\alpha, \zeta, \Delta)$ over $F$ such
that the following properties hold: $$
\left\{ \aligned \text{(i)}& \ \
\text{$R(\alpha, \zeta, \Delta)_E=
R_E(\alpha, \zeta, \Delta_E)$ for each field
extension $E/F$;}\\ \text{(ii)}& \ \
\text{$R(\alpha,
\zeta, \Delta)$ is a  simple algebra;}\\
\text{(iii)}& \ \ \text{$R$ is the
subalgebra of
$R(\alpha, \zeta, \Delta)$;}\\
\text{(iv)}& \ \  \A(R(\alpha,
\zeta, \Delta))    \text{ stabilizes $R$;}\\
\text{(v)}& \ \ \A(R(\alpha,
\zeta, \Delta)) \rightarrow \A (R),
 \sigma \mapsto \sigma|_{R} , \text{  is the
isomorphism;}\\ \text{(vi)}& \ \
\text{if } F=k \hbox{ and } \A (R_K) \hbox{ is the
$k$-group, then } \A(R(\alpha, \zeta,
\Delta)_K)  \\\vspace{-2mm}  & \hskip
2.5mm \text{is the $k$-group and } \A
(R(\alpha, \zeta,
\Delta)_K)\overset\text{(i)}\to
=\A(R_K(\alpha, \zeta,
\Delta_K))\rightarrow 
\\
 & \hskip 2.5mm \rightarrow \A(R_K), \
\sigma\mapsto \sigma|_{R_K}, \hbox{ is the
$k$-isomorphism.}
\endaligned
 \right. \tag4.2$$

\demo{{R}emark}  It follows from the
properties (ii) and (i) that the $F$-algebra
$R(\alpha, \zeta, \Delta)$ is absolutely
simple, i.e., $R(\alpha, \zeta, \Delta)_E$
is simple for each field extension $E/F$.
\enddemo

 By definition, $R$
and $Z$ are the subalgebras of $R(\alpha,
\zeta, \Delta)$, the sum of their underlying
vector spaces is direct, and there is an
element $e\in R(\alpha, \zeta, \Delta)$ such
that
$$ \mathop{\rm vect}(R(\alpha, \zeta,
\Delta))=\langle e\rangle \oplus
\mathop{\rm vect}(Z)\oplus
\mathop{\rm vect}(R) \tag4.3 $$
and the following conditions hold:
\medbreak  (R1)  $e$ is the left identity
of $R(\alpha, \zeta, \Delta)$; 
\vglue6pt \hangindent=49pt\hangafter=1 (R2) 
$\mathop{\rm vect}(Z)$ and
$\mathop{\rm vect}(R)$ in (4.3) are
respectively the eigenspaces with
eigenvalues $\zeta $ and $\alpha $ of the
operator of right multiplication of\break
$R(\alpha, \zeta, \Delta)$ by~$e$;

\vglue6pt (R3) for all $a\in R$ and $z\in Z$,
their products in $R(\alpha, \zeta, \Delta)$
are given by

\vglue-12pt

$$ \ az=0, \quad za=\Delta(z, a)e. \tag4.4
$$
\vglue4pt

Properties (4.2)(i) and (4.2)(iii)
immediately follow from this definition. Let
us show that (4.2)(ii) holds.

\nonumproclaim{Proposition 6}   The
algebra $R(\alpha, \zeta, \Delta)$ is
simple. 
\endproclaim

\demo{Proof} Let $I$ be a nonzero ideal of
$R(\alpha, \zeta, \Delta)$. Take an element
$x\in I$, $x\neq 0$. By (4.3) we have
$x=\gamma e + x^{}_Z + x^{}_R$ for some
$\gamma\in F$, $x^{}_Z\in Z$, $x^{}_R\in R$.
From (R1), (R2) we deduce that $$ I\ni
xe=\gamma ee + x^{}_Ze +x^{}_Re= \gamma e +
\zeta x^{}_Z+ \alpha x^{}_R. \tag4.5 $$

Fix an element $z\in Z$, $z\neq 0$. As
$\Delta$ is nondegenerate, there is an
element $a\in R$ such that $$ \Delta(z,
a)=1. \tag4.6 $$
As $R$ and $Z$ are the subalgebras
of $R(\alpha, \zeta, \Delta)$, formulas
(4.5), (4.1), (4.4) and condition (R1) imply
that $$ I\ni (xe)z= \gamma ez + \zeta
x^{}_Zz + \alpha x^{}_R z= \gamma z. \tag4.7
$$
From (4.7), (4.4), (4.6) we
deduce that $$ I\ni ((xe)z)a = \gamma
za = \gamma e. \tag4.8 $$
As $I$ is the ideal, (4.8)
and (R1) give $I=R(\alpha, \zeta,
\Delta)$ whenever $\gamma\neq 0$.

Consider the remaining case where $\gamma
=0$, i.e., $x=x^{}_Z+x^{}_R$. As $x\neq 0$,
either $x^{}_Z$ or $x^{}_R\neq 0$. If
$x^{}_Z\neq 0$ (resp., $x^{}_R\neq 0$), then
by nondegeneracy of $\Delta$ there is $a'\in
R$ (resp., $z'\in Z$) such that
$\Delta(x^{}_Z, a')=1$ (resp., $\Delta(z',
x^{}_R)=1$). Then since $R$ and $Z$ are the
subalgebras of $R(\alpha, \zeta, \Delta)$,
we deduce from (4.4), (4.1) that $I\ni xa' =
x^{}_Za' + x^{}_Ra'=e+a''$ for some $a''\in
R$ (resp., $I\ni z'x=z'x^{}_Z +z'x^{}_R=e$).
Thereby we have returned   to the case $\gamma
\neq 0$.   
\enddemo

The following statement immediately implies
(4.2)(iv) and (4.2)(v).

\nonumproclaim{Proposition 7} With respect to
decomposition $(4.3)$ and the
pairing~$\Delta${\rm ,} $$ \A (R(\alpha, \zeta,
\Delta))=\{ {\rm id}_{\langle e
\rangle}\!\oplus\,(g^*)^{-1}\oplus g \mid
g\in \A (R)\}. \tag4.9 $$
\endproclaim
\pagebreak

\demo{Proof} As $R$ and $Z$ are the
subalgebras of $R(\alpha, \zeta,
\Delta)$, formula (4.3),\break  Lemma 1 and
conditions (R1), (R2) imply that each
element $\sigma\in\A (R)$ has the form
${\rm id}_{\langle e
\rangle}\oplus\ h\oplus g$, $h\in
\GL(Z)$, $g\in \A (R)$. As in the
proof of Proposition 4 we obtain
$h=(g^*)^{-1}$. Thus the left-hand side
of equality (4.9) is contained in its
right-hand side.

To prove the inverse inclusion take an
element
$\varepsilon={\rm id}_{\langle e
\rangle}\oplus\ (t^*)^{-1}\oplus t$, where
$t\in \A (R)$. We have to show that $$
\varepsilon(xy)=\varepsilon(x)\varepsilon(y),\
\ x, y\in R(\alpha, \zeta, \Delta). \tag4.10
$$

If $x, y\in Z$ or $x, y\in R$, then, since
$R$ and $Z$ are the subalgebras of
$R(\alpha, \zeta, \Delta)$, equality (4.10)
follows from (4.1). If $x=e$ (resp., $y=e$),
then (4.10) follows from (R1) (resp., (R2)).
Further, (4.4) readily implies (4.10) for
$x\in R$ and $y\in Z$. Finally, if $x\in Z$
and $y\in R$, then (4.10) follows from (4.4)
by the arguments similar to those at the end
of the proof of Proposition 4; the details are
left to the reader.  
\enddemo

\nonumproclaim{{C}orollary}  Property {\rm
(4.2)(vi)} holds.  
\endproclaim

\demo{Proof} Similar to that of
the corollary of Proposition 4.  
\enddemo

 {\it Remark.} A slight
modification of the arguments makes it
possible to drop the condition
$\zeta\neq 0$, and thereby to replace
the condition $|F|\geqslant 4$ by
$|F|\geqslant 3$. That is, put $\zeta=0$
in the definition of $R(\alpha, \zeta,
\Delta)$. Then Proposition 6 holds with
the same proof. Proposition 7 holds as
well but as $\zeta =0$, Lemma 1 is not
applicable, so the proof of Proposition
7 should be modified. This can be done
as follows (with the same notation).

We have $\sigma(e)=\lambda e + e^{}_Z +
e^{}_R$ for some $\lambda\in F$, $e^{}_Z\in
Z$, $e^{}_R\in R$. As $\sigma(e)$ is the
left identity, $e=\sigma(e)e=\lambda ee +
e^{}_Ze + e^{}_Re= \lambda e + \alpha
e^{}_R$. Hence $\lambda =1$, $e^{}_R=0$. If
$e^{}_Z\neq 0$, nondegeneracy of $\Delta $
implies that $\Delta(e^{}_Z, a)=1$ for some
$a\in R$. Therefore
$a=\sigma(e)a=(e+e^{}_Z)a=a+ \Delta(e^{}_Z,
a)e=a + e$ which is impossible as $e\neq 0$.
Thus $\sigma(e)=e$. This and (R2) imply that
$\sigma(Z)=Z$, $\sigma(R)=R$. The rest of
the proof   remains unchanged.

\section{Proofs of theorems}
 
{\it Proof of Theorem $1$}. Let $U$,
$b$, $r$, $L$ and $S$ be as in the
formulation of Proposition~5. We may
(and shall) assume that   $r>1$. The
vector spaces $U$, $L$ and $S$ being
defined over $k$, let $U_0$, $L_0$ and
$S_0$ be the corresponding
$k$-structures. Put $s:=\dim L$,
$n:=\dim U$.

Assume that the number of elements in
$k$ satisfies inequality (2.14) for
$F=k$. Then we may (and shall) fix some
$\g$, $\de$, $\Phi$, $\alpha$, $\zeta$,
$\Delta$ and consider the $k$-algebra
$A:=R(\alpha, \zeta, \Delta)$, where
$R:=D(L_0, U_0, S_0, \g, \de, \Phi)$
(see Sections 2 and 4).

It follows from Proposition 5,
Corollary of Proposition 4 and property
(2.18) that $\A(R_K)$ is the $k$-group
$k$-isomorphic to $G$. Now the claim
follows from Proposition 6 and
Corollary of Proposition 7.  
\hfill\qed

\demo{Proof of Theorem $2$} By Theorem 1,
there is a finite dimensional algebra $A$
over $k$ such that $\A(A_K)$ is the
$k$-group $k$-isomorphic to $G$. Put
$V:=\mathop{\rm vect}(A)$. Then
$V^*_K\otimes V^*_{K}\otimes V^{}_{K}$ is
the variety of all $K$-algebra structures
(i.e., multiplications) on $V_K$ (to $\sum
f\otimes l\otimes v\in V^*_K\otimes
V^*_{K}\otimes V^{}_{K}$ corresponds the
multiplication defined by $xy=\sum f(x)l(y)
v, \ x, y\in V_K$). Two $K$-algebra
structures correspond to isomorphic algebras
if and only if their $\GL(V_K)$-orbits
coincide. In particular the
$\GL(V_K)$-stabilizer of a tensor $t\in
V^*_K\otimes V^*_{K}\otimes V^{}_{K}$ is the
full automorphism group of the algebra
corresponding to $t$, cf.\hskip
1mm\cite{Se}. Therefore the tensor
corresponding to multiplication in $A$ is
the one we are after.  
\enddemo

\section{Constructive proof of Corollary
2 of Theorem 1 }
 
Using the fact  that regular representation of a
finite abstract group $G$ yields its
faithful representation by permutation
matrices, one immediately deduces
Corollary~2 from Theorem 1. Since our
proof of Theorem 1 is nonconstructive,
this proof of Corollary 2 is
nonconstructive as well. Here we give
another, constructive proof of this
corollary. Combined with our proof of
Theorem 2, it yields a constructive
realization of $G$ as the
$\GL(V_K)$-stabilizer of a $k$-rational
tensor in $V^*_K\!\otimes
\!V^*_{K}\!\otimes \!V^{}_{K}$ for some
finite dimensional vector space $V$
over $k$. Our constructive proof works
if $k$ contains sufficiently many
elements (for instance, if $k$ is
infinite).

By Lemma 2, we may (and shall) assume
that $G$ is nontrivial. For an
appropriate $n>1$ fix an embedding
$\iota\!: G\hookrightarrow \s$ and
identify $G$ with the subgroup $\iota
(G)$ of $\s$. Let $E_n$ be the
$n$-dimensional split \'etale algebra
over $k$, i.e., the  direct sum of $n$
copies of the field $k$. Put
$V:=\mathop{\rm vect}(E_n)$ and
denote by $e_i$ the 1 of the $i^{\rm th}$ direct
summand of $E_n$. Consider the natural
action of $\s$ on $E_n$ and $V$
given~by $$ \sigma\cd
e_i=e_{\sigma(i)}, \ \ \sigma \in \s,
1\leqslant i \leqslant n. \tag6.1 $$
 As $\s$ acts faithfully, we
may (and shall) identify $\s$ with the
subgroup of $\GL(V)$. Then it is easily
seen that $\A (E_n)=\s$.

If $k$ contains sufficiently many
elements, one can find a sequence of
nonzero elements
$\hskip-8pt\mathbold{\lambda}:=(\lambda_1,
\ldots, \lambda_n)\in k^n$ such that
 $$
\lambda_i/\lambda_j\neq\lambda_s/\lambda_t\
\text{ for all } 1\leqslant i, j, s,
t\leqslant n, \ i\neq j, s\neq t, (i,
j)\neq (s, t). \tag6.2 $$
  Fix such a sequence.

\nonumproclaim{Proposition 8}   Put $
\textstyle f:=\prod_{\sigma\in
G}\sigma\cd(\lambda_1e_1+{\cdots}+\lambda_ne_n)
\in \Sy^{|G|}(V). $ Then $G=(\s)_{\langle
f\rangle}.$ 
\endproclaim

\demo{Proof} The definition of $f$
clearly implies that $G\subseteq
(\s)_{\langle f\rangle}$. To prove the
inverse inclusion, take an element
$\delta\in (\s)_{\langle f\rangle}$. As
$\lambda_1e_1+{\cdots}+\lambda_ne_n\in\Sy(V)$
divides $f$, and $\delta\in
\A(\Sy(V))$, we deduce that
$\delta\cd(\lambda_1e_1+{\cdots}+\lambda_ne_n)$
divides $\delta\cd f$. As the algebra
$\Sy(V)$ is factorial and
$\lambda_1e_1+{\cdots}+\lambda_ne_n$ is
its prime element, this, plus the definition
of $f$ and the inclusion $\delta\cd
f\in \langle f\rangle$ yield $$
\delta\cd(\lambda_1e_1+{\cdots}
+\lambda_ne_n)= \alpha(\sigma \cd
(\lambda_1e_1+{\cdots}+\lambda_ne_n)) \
\ \text{ for some } \sigma\in G, \
\alpha\in k. \tag6.3 $$

We claim that $\delta=\sigma$. If not,
there are two indices $i_0$ and $j_0$,
$i_0\neq j_0$, such that $$
\delta^{-1}(i_0)\neq \sigma^{-1}(i_0) \
\text{ and } \delta^{-1}(j_0)\neq
\sigma^{-1}(j_0). \tag6.4 $$
  From (6.3) and (6.1) we obtain
$\lambda_{\delta^{-1}(i)}
=\alpha\lambda_{\sigma^{-1}(i)}$ for each
$i=1,{\ldots}\,,n$. Hence
$\lambda_{\delta^{-1}(i_0)}/
\lambda_{\delta^{-1}(j_0)}=
\lambda_{\sigma^{-1}(i_0)}/
\lambda_{\sigma^{-1}(j_0)}$. By (6.2), this
implies $\delta^{-1}(i_0)=\sigma^{-1}(i_0)$,
$\delta^{-1}(j_0)=\sigma^{-1}(j_0)$,
contrary to (6.4).  
\enddemo

 {\it Remark.} It follows from
the definition of $f$ that $G\subseteq
(\s)\!_f$. Hence
$(\s)\!_f=(\s)_{\langle f\rangle}=G$.
\vglue9pt

Assuming that $k$ contains sufficiently many
elements, fix a sequence
$\hskip-8pt\mathbold{\mu}:=(\mu_1,\dots,\mu_{|G|+1})\in
k^{|G|+1}$ where $\mu_i\in k\setminus\{0,
1\}$, $\mu_i\neq \mu_j$ for $i\neq j$, and
define the algebra $E(k, G, \iota,
\hskip-8pt\mathbold{\lambda}, \hskip-8pt\mathbold{\mu})$ as
follows.

Put $S:=\langle f\rangle\subset
\mathop{\rm Sym}^{|G|}(V)$. First, $A(V,
S)$ and $E_n$ are the subalgebras of $E(k,
G, \iota, \hskip-8pt\mathbold{\lambda},
\hskip-8pt\mathbold{\mu})$ and the sum of their
underlying vector spaces is direct. So
$\mathop{\rm vect}(E(k, G, \iota,
\hskip-8pt\mathbold{\lambda}, \hskip-8pt\mathbold{\mu}))$
contains two distinguished copies of $V$:
the copy $V^{}_A$ corresponds to the summand
$\Sy^{1}(V)$ in (2.3), and the copy $V^{}_E$
to the subalgeb\-ra~$E_n$.

Let $(\ {,}\ )\!: V^{}_A\times V^{}_E
\rightarrow k$ be the nondegenerate bilinear
pairing defined by the condition
$$\gathered
(e_i,e_j)=\cases 1\ \text{
if } i=j,\\ 0 \ \text{ if } i\neq j.
\endcases
\endgathered
\tag6.5 $$

Second, there is an element $e\in E(k, G,
\iota, \hskip-8pt\mathbold{\lambda}, \hskip-8pt\mathbold{\mu})$
such that $$ \mathop{\rm vect}(E(k, G,
\iota, \hskip-8pt\mathbold{\lambda}, \hskip-8pt\mathbold{\mu}))
= \langle e \rangle\oplus
\mathop{\rm vect}(A(V, S))\oplus
\mathop{\rm vect}(E_n)\tag6.6
$$
and  the following conditions
hold:

\medbreak \item{(E1)} $e$ is the left identity
of $E(k, G, \iota, \hskip-8pt\mathbold{\lambda},
\hskip-8pt\mathbold{\mu})$; 
\pagebreak \item{(E2)}
$\mathop{\rm vect}(E_n)$ in (6.6) and each
summand in decomposition (2.3) of the\break
summand $\mathop{\rm vect}(A(V, S))$ in
(6.6) are respectively the eigenspaces with
eigenvalues $\mu_1, \ldots, \mu_{|G|+1}$ of
the ope\-ra\-tor of right multiplication of\break
$E(k, G, \iota, \hskip-8pt\mathbold{\lambda},
\hskip-8pt\mathbold{\mu})$ by~$e$; 
\vglue8pt\item{(E3)} if $x$
is an element of a direct summand in (2.3),
and $y\in V^{}_E$, then $$ xy=yx= \cases (x,
y)e
& \text{ if } x\in V^{}_A,\\
\hskip 4mm 0 & \text { otherwise,}
\endcases $$
\vglue8pt

It readily follows from the definition
that for each field extension $F/k$ we
have  $$ E(k, G, \iota,
\hskip-8pt\mathbold{\lambda}, \hskip-8pt\mathbold{\mu} )_F=
E(F, G, \iota, \hskip-8pt\mathbold{\lambda},
\hskip-8pt\mathbold{\mu}).\tag6.7 $$

\phantom{great}

\nonumproclaim{Proposition 9}   With
respect to decomposition $(6.6)${\rm ,}  $$
\A (E(k, G, \iota, \hskip-8pt\mathbold{\lambda},
\hskip-8pt\mathbold{\mu} ))=
\{{\rm id}_{\langle
e\rangle}\!
{\oplus}\,g\,{\oplus}\,g\mid g\in G\}.
\tag 6.8 $$
\endproclaim

\demo{Proof} Lemma 1, Proposition 1 and
conditions (E2), (E3), (E1) yield that
with respect to decomposition (6.6)
each element of $\A (E(k, G, \iota,
\hskip-8pt\mathbold{\lambda}, \hskip-8pt\mathbold{\mu}))$
has the form
${\rm id}_{\langle
e\rangle}\!{\oplus}g{\oplus}h$ for some
$g\!\in\! \GL(V)_S$, $h\!\in\! \s$.
As in the proof of Pro\-po\-sition~4
we obtain $h=(g^*)^{-1}$. On the other
hand, (6.1), (6.5) imply that $$
(\sigma^*)^{-1}=\sigma \ \text{ for
each } \sigma \in \s. \tag6.9 $$

\noindent Now Proposition 8 and (6.9)
yield that $g\in \GL(V)_S \cap
\s=(\s)_S= G$. Thus the left-hand side
of equality (6.8) is contained in its
right-hand side. The inverse inclusion
is verified as at the end of proof of
Proposion 4 (with (6.9) taken into
account); we leave the details to the
reader.  
\enddemo

From (6.7) and (6.8) we immediately
deduce the following.

\nonumproclaim{{C}orollary  1}   The groups $G$
and $\A(E(k, G, \iota, \hskip-8pt\mathbold{\lambda},
\hskip-8pt\mathbold{\mu})_F)$ are isomorphic for each
field extension $F/k$. 
\endproclaim

In turn, this implies the constructive proof
of Corollary 2 of Theorem 1. Namely, assume
that $k$ contains sufficiently many
elements, fix some $\hskip-8pt\mathbold{\lambda}$,
$\hskip-8pt\mathbold{\mu}$, $\iota$, $\alpha$,
$\zeta$, $\Delta$ and consider the
$k$-algebra $A:=R(\alpha, \zeta, \Delta)$
(see Section 4), where $R:=E(k, G, \iota,
\hskip-8pt\mathbold{\lambda}, \hskip-8pt\mathbold{\mu})$. Then
properties (4.2)(i),(ii),(v) and Corollary 1
yield the following.

\nonumproclaim{{C}orollary  2}   The finite
dimensional $k$\/{\rm -}\/algebra $A$ is simple and
$\A(A_F)$ is isomorphic to $G$ for each
field extension $F/k$. 
\endproclaim

\section{Appendix}

Realization of algebraic groups as
normalizers and stabilizers is crucial
for this paper: our proof of Theorem 1
is based on realization of algebraic
groups as normalizers of some linear
subspaces; Theorem 2 concerns
realization of algebraic groups as
stabilizers of some very specific
tensors. This appendix contains further
results on this topic. In particular it
yields a refinement of Proposition~5.

\nonumproclaim{Proposition 10}  Let $G$
be an algebraic $k$\/{\rm -}\/group. There are an
integer  $n>0$ and a closed
$k$\/{\rm -}\/embedding $G\hookrightarrow
R:=\GL(1, K)\times \GL(n, K)$ such that
for each closed $k$\/{\rm -}\/embedding of $R$
in an algebraic $k$\/{\rm -}\/group $Q$ the group
$G$ is the stabilizer of a $k$\/{\rm -}\/rational
element of a finite dimensional
$Q$\/{\rm -}\/module defined over~$k$.
\endproclaim

{\it Proof}. As $G$ is algebraic, we may
(and shall) consider it is as a closed
$k$-subgroup of $\GL(n, K)$ for some $n$. By
Chevalley's theorem, cf.\ [H, 11.2,
34.1], [Sp, 5.5.3], there are a finite
dimensional $\GL(n, K )$-module $U$ and a
one-dimensional linear subspace $S$ of $U$,
both defined over $k$, such that\break $G=\GL(n, K
)\!_S$. We have $g \cd s=\chi(g)s,\ g\in G,
s\in S$, for some character $\chi\! :
G\rightarrow \GL(1, K)$ defined over $k$.
Consider the reductive $k$-group $R:=\GL(1,
K) \times \GL(n, K)$ and define the
$R$-module structure on $U$ by $(\lambda,
g)\cd u:=\lambda (g\cd u), \ \lambda\in
\GL(1, K), g\in\!\GL(n, K), u\!\in U$. The
$R$-module $U$ is defined over $k$. For
eve\-ry $s\in S$, $s\neq 0$, we have $R_s =
\{(\chi (g^{-1}), g) \mid g\in G\}$. Hence
$G\rightarrow R$, $g\mapsto (\chi (g^{-1}),
g)$, is the closed $k$-embedding whose image
is $R_s$.

Let $R$ be the closed $k$-subgroup of
some $Q$. By Hilbert's theorem,
cf.\ \cite{MF}, \cite{PV}, as $R$ is
reductive, the homogeneous space $Q/R$
is affine. Hence $R$ is an  observable
subgroup of $Q$, cf.\ \cite{BHM}. As
$R_s$ is  an $R$-stabilizer of a vector in
an $R$-module, $R_s$ is an observable
$k$-subgroup of $R$; cf.\ \cite{BHM}.
This implies that $R_s$ is an observable
$k$-subgroup of $Q$, which in turn
implies existence of a finite
dimensional $Q$-module $M$ defined over
$k$ such that $R_s$ is a stabilizer of a
$k$-rational element of $M$;
cf.\ \cite{BHM}, [PV,\,1.2, 3.7].
\hfill\qed\vglue9pt

The following statement was used in the
first version of our proof of Theorem~1
found in the summer of 2001. Combined with
Proposition 10, it yields that in
Proposition 5 one may take
 $\dim S=1$ and $\dim
L\geqslant 1$. We believe that it is of
interest in its own right and might be
useful in other situations.

Fix a nonzero finite dimensional vector
space $U$ over $K$ endowed with a
$k$-structure $U_0$.

\nonumproclaim{Proposition 11}   Let $M$
be a finite dimensional
$\SL(U)$\/{\rm -}\/module defined over $k$ and
let $L$ be a nonzero trivial
$\SL(U)$\/{\rm -}\/module defined over $k$. Then 
 
\item{\rm (i)} 
$\T(U)_+$ contains a submodule defined
over $k$ and $k$\/{\rm -}\/isomorphic to $M${\rm ;}
\pagebreak
\item{\rm (ii)} there
is an integer $b>0$ such that $(L\oplus
U)^{\otimes m}$ contains a submodule
defined over~$k$ and $k$\/{\rm -}\/isomorphic to
$M$  for each $m\geqslant b$. 

\endproclaim

{\it Proof}. (i) Let $\Cal M$ be the
class of finite dimensional submodules
of $\T(U)_+$.

Let  $u_1,{\ldots}\,,u_n$ be a
$k$-rational basis of $U$. For every
$d=1, \ldots, n$  put
$t_d:=\sum_{\sigma\in\frak{S}_d
}\mathop{\rm sgn}(\sigma)
u_{\sigma(1)}\otimes\cdots\otimes u_
{\sigma(d)}$. This is a $k$-rational
skew-symmetric element of $U^{\otimes
d}$, fixed by $\SL(U)$ for $d=n$. Hence
for all $m, l
>0$ the module $U{}^{\otimes (m+ln)}$
contains the submodule $t^{\otimes
l}_n\otimes U{}^{\otimes m}$ defined
over $k$ and isomorphic to $
U{}^{\otimes m}$ over $k$. This implies
that if $M_1, M_2\!\in\! \Cal M$
(resp., $M_1, M_2\in \Cal M$ and are
defined over $k$), then
$M_1\nomathbreak\oplus\nomathbreak M_2$
and $M_1\otimes M_2$ are isomorphic
(resp., isomorphic over $k$) to the
modules from $\Cal M$ (resp., defined
over $k$). The class $\Cal M$ is also
closed under taking submodules. In
particular it is closed under taking
direct summands. On the other hand,
 $\Cal M$ contains all `fundamental'
 $\SL(U)$-modules: the $d^{\rm th}$ one is
 the minimal submodule of $\T(U)_+$ containing
$t_d$. If $\ch k=0$, using complete
reducibility and `highest weight
theory', we immediately deduce from
these facts that every
finite dimensional $\SL(
U)$-module is isomorphic to a module
from $\Cal M$. We claim that this is
true for $\ch k>0$ as well. To prove
this, and then to complete the proof of
(i), we use the arguments communicated
to us by W.\,van\,der Kallen,
[vdK2] . We use the standard
notation,~\cite{J}.

The above arguments show that, for
$p:=\ch k>0$, each tilting module is
isomorphic to a module from $\Cal M$
(see the necessary information on
titlting modules at the end of this
section).

\nonumproclaim{Lemma 3  {\rm
(Cf.\,\cite{Don})}}  For $m$
sufficiently large{\rm ,}
$M\otimes\St_m\otimes\St_m$ is tilting.
\endproclaim

{\it Proof}. As $\St_m$ is self-dual,
it suffices to show that $M\otimes\St_m
\otimes\St_m$ has good filtration for
$m$ large. In fact we will show that
$M\otimes \St_m$ has good filtration
for $m$ large. Take $m$ so large that
for each weight $\lambda$ of $M$ the
weight $\lambda+(p^m-1)\rho$ is
dominant. Then $M\otimes (p^m-1)\rho$
is what is called in \cite{Pol} a
module with excellent filtration for
the Borel group $B$ whose roots are
negative; cf.\ [vdK1]. Indeed
it has a filtration with each  layer
one-dimensional of dominant weight.
Therefore, by Kempf vanishing,
$M\otimes \St_m =
\mathop{\rm ind}_B^G (M\otimes
(p^m-1)\rho)$ has good filtration.  
\hfill\qed\vglue9pt

As $\St_m\otimes \St_m$ contains a
trivial one-dimensional submodule,
Lemma 3 shows that $M$ can be embedded
in a tilting module, whence the claim.

Now we take into account the field of
definition $k$. According to what we
proved, $M$ is isomorphic to a
submodule of some $N:=\hskip-8pt\mathbold{\oplus}_{i=1}^{r}
{U}{}^{\otimes i}$. Let $M_0\subset M$
be the $k$-structure of $M$. As
$N_0:=\hskip-8pt\mathbold{\oplus}_{i=1}^{r}U^{\otimes i}_0$
is the $k$-structure\break of $N$,
$\mathop{\rm Hom}_{\SL(U_0)}(M_0,
N_0)$ is the $k$-structure of
$\mathop{\rm Hom}_{\SL( U)}(M, N)$;
cf.\ [J, I, 2.10(7)]. So there are
$h_i\in
\mathop{\rm Hom}_{\SL(U_0)}(M_0,
N_0)$ and $\lambda_i\in K$,
$i=1,{\ldots}\,,m$, such that
$\sum_{i=1}^{m}\lambda_i({h}_i\otimes
1)\!: M\!\rightarrow\! N$ is an
injection of $\SL( U)$-modules. But
this implies that
$\hskip-8pt\mathbold{\oplus}_{i=1}^{m}({h}_i\otimes\nomathbreak
1)\!: M\rightarrow N^{\oplus m}$ is the
injection of $\SL( U)$-modules defined
over $k$. Now we remark that, according
to what is already proved, as $N$ is a
submodule of $\T( U)_+$ defined over
$k$, there is an injection $N^{\oplus
m}\rightarrow \T( U)_+$ defined over
$k$. This completes the proof of (i).

\vglue4pt

(ii) As $M$ is finite dimensional, by
(i) there is  an integer $b$ such that
$M$ is isomorphic over $k$ to a
submodule of $\bigoplus_{i=1}^{b}
U{}^{\otimes i}$ defined over $k$. Take
a nonzero $k$-rational element $l\in
L$. Then, for each $i$, $1\leqslant
i\leqslant b$, the $\SL( U)$-module $
U{}^{\otimes i}$ defined over $k$ is
isomorphic over $k$ to the submodule
$l^{\otimes (m-i)}\otimes
 U{}^{\otimes i}$ of $(L\oplus
{U}){}^{\otimes m}$ defined over $k$
(here $l^{\otimes 0}\otimes
U{}^{\otimes m}$ stands for $
U{}^{\otimes m}$), whence the claim. \phantom{greatmusic}\hfill\qed
\vfil

{\elevensc Corollary}.  {\it Let $G$ be an
algebraic $k$\/{\rm -}\/group. Then there is a
finite dimensional vector space $U$
defined over $k${\rm ,} a closed
$k$\/{\rm -}\/embedding $G\hookrightarrow \SL(
U)$ and an integer $b>0$ such that for
each integer $m\geqslant b$ and nonzero
trivial $\SL(U)$\/{\rm -}\/module $L$ defined
over $k$ the group $G$ is the
$\SL(U)$\/{\rm -}\/stabilizer of a $k$\/{\rm -}\/rational
tensor in $(L\oplus
 U){}^{\otimes m}$.
}\vfil

{\it Proof}. As each algebraic
$k$-group is a closed $k$-subgroup of
$\SL( U)$ for some vector space $U$
defined over $k$, the claim follows
from Propositions 10 and 11. 
\phantom{hi}\hfill\qed\vfil

{\it Tilting modules.}
For the reader's convenience,  we
collect here the basic definitions and
(nontrivial) facts about tilting
modules used in the proof of
Proposition 11; cf.\  \cite{Don} and
the references therein.

Let $G$ be a reductive  connected linear 
algebraic group, $T$ a maximal torus
and $B\supseteq T$ a Borel subgroup of
$G$. Let
${\rm X}={\rm X}(T)$
be the character group of~$T$. We fix
in ${\rm X}$ the system of
simple roots of $G$ which makes $B$ the
negative Borel. For $\lambda\in
{\rm X}$ we denote by
$K_{\lambda}$ the one-dimensional
$B$-module on which $T$ acts with
weight $\lambda$. Let $\nabla(\lambda)$
be the induced $G$-module
$\mathop{\rm Ind}^{G}_{B}(K_{\lambda})$
(i.e., the $G$-module of global
sections of the homogeneous line bundle
over $G/B$ with the fiber $K_{\lambda}$
over $B$). Then $\nabla(\lambda)$ is
finite dimensional and is nonzero if
and only if $\lambda $ belongs to the
monoid ${\rm X}^{+}$ of
dominant weights; cf.\ \cite{J}.

An ascending filtration of a $G$-module
is called {\it good} if each successive
quotient is either zero or isomorphic
to $\nabla(\lambda)$ for some
$\lambda\in {\rm X}^{+}$. Let
$\Cal T$ be the class of finite
dimensional $G$-modules $V$ such that
both $V$ and its dual $V^*$
have good filtration. Then one can prove the
following.
\vfil
(i) A direct sum and a tensor product
of modules in $\Cal T$ belong to $\Cal
T$ and also, a direct summand of a
module in $\Cal T$ belongs to $\Cal T$.

\vfil
(ii) For each $\lambda\in
{\rm X}^{+}$ there is an
indecomposable (into a direct sum) module
$M(\lambda)\in \Cal T$ which has unique
highest weight $\lambda$; furthermore,
$\lambda $ occurs with multiplicity 1 as a
weight of $M(\lambda)$, and the modules
$M(\lambda)$, $\lambda\in
{\rm X}^{+}$, form a complete set
of inequivalent indecomposable modules in
$\Cal T$; cf.\ \cite{Don}.
\vfil

The mo\-dule $M(\lambda)$ ($\lambda\in
{\rm X}^{+}$) is called the {\it
tilting} module of highest weight~$\lambda
$.  \eject

 \AuthorRefNames [KMRT]
\references 
 
 \vglue12pt
[BHM]
\name{A.\ Bialynicki-Birula, G.\ Hochshild}, and \name{G.\ D.\  Mostow},
  Extensions of representations of algebraic linear groups,  {\it Amer.\ J.\ Math.\/}
(1963)
{\bf 65}, 131--144.

[Bor]
 \name{A.\ Borel},
{\it Linear Algebraic Groups}, 
2nd ed., {\it Graduate Texts in Math.\/} {\bf 126},\break
  Springer-Verlag, New York, 1991.

[Don]  \name{S.\ Donkin},  On
tilting modules for algebraic groups,
{\it Math.\ Z.} {\bf 212} (1993), 39--60.

[Dix]
 \name{J.\ Dixmier},
  Certaines alg\graveaccent ebres non
associatives simples d\'efinies par la
trasvections du formes binaires,
{\it J.\ fur die reine und angew.\ Math.\/}
{\bf 346} (1984),
 110--128.

[DG]
 \name{M.\ Dugas} and \name{R.\  G\"obel},
 Automorphism groups of fields, {\it Manuscripta Math.\/}
{\bf 85} (1994), 227--242.

[F]
\name{M.\ Fried},
 A note on automorphism groups of
algebraic number fields, {\it Proc.\ Amer.\ Math.\ Soc.\/}
{\bf 80}
(1980),  386--388.

[Ge]
 \name{W.\  D.\ Geyer},
 Jede endliche Gruppe ist Automorphismengruppe
einer endlichen Erweiterung $K/\bold Q$,
{\it Arch Math.\/}
{\bf 41}
(1983), 139--142.

[Gr]
 \name{R.\ L.\  Griess, Jr.},
  The friendly giant,
{\it Invent.\ Math.\/} {\bf 69}
(1982), 1--102.

[H]
 \name{J.\ E.\  Humphreys},
{\it Linear Algebraic Groups},  3rd ed.,
{\it Graduate Texts in Math.\/} {\bf 21},
  Springer-Verlag, New York, 1987.

[Ilt]  \name{A.\ V.\   Iltyakov},
 Trace polynomials and invariant
theory, {\it  Geom.\  Dedicata} {\bf 58}
(1995), 327--333.

[J]
 \name{J.\ C.\  Jantzen},
{\it  Representations of Algebraic Groups{\rm ,}
 Pure and Applied Mathematics} {\bf 131},
 Academic Press, Inc.,  Boston, 1987.

[K1]  \name{V.\ Kac}   Some
remarks on nilpotent orbits, {\it  J.\
Algebra} {\bf 64 } (1980), 190--213.

[K2] \bibline,   Letter of
March  2, 2002.

[KMRT]
  \name{M.-A.\ Knus, A.\
Merkurjev, M.\ Rost}, and \name{J.-P.\ Tignol},
{\it The Book of Involutions}, 
{\it Colloq.\  Publ.\ Amer.\ Math.\ Soc.\/} {\bf 
44},  Providence, RI,  1998.

[MFK]
\name{ D.\ Mumford, J.\
Fogarty}, and \name{F.\  Kirwan},
{\it Geometric
Invariant Theory}, {\rm 3rd enl.\ ed.},
{\it Ergeb.\  Math.\ 
Grenzgeb.\/} 
{\bf 34},
Springer-Verlag, New York,
1994.

[Pol]
 \name{P.\ Polo},
  Vari\'et\'es de Schubert et excellentes filtrations, in
  {\it Orbites Unipotentes et Representations}.\  III,
{\it Ast{\rm \'{\it e}}risque}
{\bf 173--174}
 (1989),
 281--311.

[Pop]  \name{V.\ L.\  Popov},  
An analogue of M.\ Artin's conjecture on
invariants for non-associa\-ti\-ve
algebras, {\it Amer.\ Math.\ Soc.\ Transl.\/} {\bf 169}
(1995), 121--143.

[PV]
 \name{V.\ L.\  Popov} and \name{E.\ B.\  Vinberg},
 Invariant theory, 
in {\it Algebraic Geometry} IV, {\it Encycl.\ Math.\ Sci.\/} {\bf 55}, 123--284,
 Springer-Verlag, New York, 1994.
 
[Re]
 \name{Z.\ Reichstein},
 On the notion of essential dimension for
algebraic groups,  {\it Transformation Groups}
{\bf 5}
(2000), 265--304.

[Sc]
 \name{R.\ D.\  Schafer},
{\it Introduction to Nonassociative Algebras},
 Academic Press, New York, 1966.

[Se]
 \name{J-P.\ Serre},
{\it Galois Cohomology},
 Springer-Verlag, New York,
  1997.

[Sp]  \name{T.\ A.\  Springer},  {\it Linear Algebraic Groups},  Second Edition,
 {\it Progr.\  Math.\/} {\bf 9},\break
Birkh\"auser, Boston ,  1998.

[SV]  \name{T.\ A.\  Springer}, and \name{F.\ D.\  Veldkamp}  
{\it Octonians{\rm ,} Jordan Algebras and Exceptional Groups}, 
Springer-Verlag, New York, 2000.

[vdK1]
  \name{W.\  van der
Kallen},
{\it  Lectures on Frobenius
Splittings and $B$-Modules}. {\it Notes by S.\
P.\  Inamdar}, 
  Published for the Tata
Institute of Fundamental Research,
Bombay, by Springer-Verlag, Berlin,   1993.

[vdK2]
 \name{W.\  van der
Kallen},
  Letters of May {\rm 25} and August {\rm 2, 2001}.

\endreferences
\bye